\documentclass{article}
\usepackage{amssymb}
\usepackage{mathrsfs}
\usepackage{CJK}
\usepackage{bbm}
\usepackage{stmaryrd}
\usepackage{amsmath}
\usepackage{cases}
\usepackage{amscd}
\usepackage{amsfonts}
\usepackage{latexsym,amssymb,amsmath,mathrsfs,amsbsy, amsthm}
\usepackage[usenames]{color}
\usepackage{xspace,colortbl}
\usepackage[
pdfstartview=FitH, CJKbookmarks=true, bookmarksnumbered=true,
bookmarksopen=true,
citecolor=blue %¸Ä±äĿ¼¡¢ÒýÓõÄÑÕÉ«
,colorlinks=true
,pdfborder=1 %×¢Ê͵ô´ËÏîÔò½»²æÒýÓÃΪ²ÊÉ«±ß¿ò
,pdfstartpage=1 %´ò¿ªÊ±ÎªµÚ2Ò³
,pdfstartview=Fit]{hyperref} %´ò¿ªpdfʱΪ fit page
\allowdisplaybreaks
\usepackage{pdfsync}

\newtheorem{thm}{Theorem}[section]
\newtheorem{lem}[thm]{Lemma}
\newtheorem{prop}[thm]{Proposition}
\newtheorem{coro}[thm]{Corollary}
\newtheorem{defn}[thm]{Definition}
\theoremstyle{definition}
\newtheorem{rem}[thm]{Remark}
\newtheorem{exmp}[thm]{Example}

\newcommand{\mz}{\mathbb{Z}}
\newcommand{\mc}{\mathbb{C}}

\numberwithin{equation}{section}
\begin{document}
\title{\bf On the Wu classes of the topological blow up %\thanks{}
 }

\author{\bf Wang Wei \\{\small School of Mathematical Sciences, Fudan University }\\
{\small Shanghai, P. R. China, 200433}}
\date{\small July, 2011}
\maketitle
\begin{abstract}
\small  We mainly discuss the Wu classes $v(\tilde{M})$ and the Steenrod operation $Sq$ of the topological blow up $\tilde{M}$. The formula of the Wu class $v(\tilde{M})$ will be given as well as the formula of the Steenrod operation $Sq$. As an application, we will use our results to describe a geometric obstruction.
\end{abstract}

\section{Introduction}
In algebraic and symplectic geometry, blow up is an important surgery to obtain new manifolds and an effective technique to resolve singularities. For the geometric and topological properties of algebraic and symplectic blow up, there are lots of references. For example, \cite{Duan},\cite{Geiges},\cite{Griff},\cite{Lascu},\cite{Macduff},\cite{Voisin}.

Topologically, these is an analogous construction. Let $M$ be a closed, smooth manifold and $N$
be a closed submanifold of $M$ with codimension $2r>2$. Suppose the normal bundle $\eta_N$ of $N$ admits a complex structure, topologically, we can blow up $M$ along $N$ to get a closed, smooth manifold $\tilde{M}$ and there exists a natural degree one map $\pi:\tilde{M}\longrightarrow M$ such that:\\
(1) $\pi: \pi^{-1}(M-N)\longrightarrow M-N$ is a diffeomorphism. \\
(2) $\pi: \pi^{-1}(N)=P(\eta_N)\longrightarrow N$ is the projection map, where $P(\eta_N)$ is the complex projectivization of the normal bundle $\eta_N$.

By the cohomolgoy property of the degree one map $\pi$,
we'll see the cohomology group $H^*(\tilde{M},\mz)$ of $\tilde{M}$ has a direct sum decomposition:
\[
H^* (\tilde{M},\mz)=\pi^* H^* (M,\mz)\oplus \text{Ker} \alpha_{\pi}^{*}
\]
where $\alpha_{\pi}^{*}$ is the Gysin map induced by $\pi$.

In section 2, after the introductions to degree one map and Gysin map in section 2.1 and the details of the construction of topological blow up in section 2.2, we will prove:
\begin{prop}
$Ker\alpha_{\pi}^{*}= \alpha_{j_P}^{*}(\oplus_{i=0}^{r-2}H^*(N,\mz){\xi}^i)$, where $\alpha_{j_P}^{*}$ is the Gysin map induced by the inclusion $j_P: P(\eta_N)\hookrightarrow \tilde{M}$ and $-\xi\in H^2(P(\eta_N),\mz)$ is the Euler class of the normal bundle of $P(\eta_N)$.
\end{prop}
As a quick application, some $L$-classes of $\tilde{M}$ will be concerned. In section 2.3, we'll see:
\begin{prop}
Suppose $\tilde{M}$ is the topological blow up of $M$ along the submanifold $N$, when $4i < min\{2r,\frac{n-1}{2}\}$,
$$L_i(\tilde{M})=\pi^{*}L_i(M)+ \alpha_{j_P}^{*}(\xi_{Q}^{-1}(L_i(S)-i_{N}^{*}(L_i (M))))$$
here $S$ is the normal bundle of $P(\eta_N)\subset \tilde{M}$ with $L_i (S)= L_i (p_1(S),p_2(S),\cdot\cdot\cdot)$ and $\sum_{i}(-1)^i p_i(S)=c(S) \overline{c(S)}$ where
\[
c(S)=(1-\xi)(\sum_{i=0}^{r}c_i(\eta_N)(1+\xi)^{r-i}) \in H^*(S,\mz)=H^*(P(\eta_N),\mz)
\]
\end{prop}
In section 3, we mainly discuss the Steenrod operation $Sq$ on $H^*(\tilde{M},{\mz}_2)$ and the Wu class $v(\tilde{M})$ of $\tilde{M}$. We will recall some facts of Wu class first, for example, Wu formula and the behavior of Wu class under degree one map. In this section, two main results will be given.

For the Steenrod operation, we have:
\begin{prop}
For any element $x=\pi^* (x_1)+ \alpha_{j_P}^{*}(x_2) \in \pi^*H^*(M,\mz_2)\oplus \alpha_{j_P}^{*}(\oplus_{i=0}^{r-2}H^*(N,\mz_2){\overline{\xi}}^i)=H^*(\tilde{M},\mz_2)$:
\[
Sq(x)=\pi^* Sq(x_1) + \alpha_{j_P}^{*}((1+\overline{\xi})Sq(x_2))
\]
\end{prop}

%By the property of the degree one map, the Wu class satisfy $v(\tilde{M})=\pi^* v(M)+ \overline{v}$, where $\overline{v}\in$ Ker$\alpha_{\pi}^{*}$.
For the Wu class, we have formula:
\begin{thm} The Wu class of the topological blow up $\tilde{M}$ has the form:
\[
v(\tilde{M})=\pi^* v(M)+ \alpha_{j_P}^{*}({\overline{\xi}}^{-1}(p^*v(N)\chi(Sq)(\sum_{i=0}^{r}w_{2i}(\eta_N)(1+\overline{\xi})^{r+1-i})-i_{N}^{*}v(M)) )
\]
\end{thm}

The formula of $v(\tilde{M})$ seems not easy to compute and in some case we only need to determine whether $v_m (\tilde{M})$ is vanished or not, so in section 3.3, we will prove:
\begin{prop}
$v_{m}(\tilde{M})=0$ if and only if\\
(1). $v_m (M)= 0$\\
(2). ${Sq}^m b + (Sq^{m-2}b) \overline{\xi}=0$ for all $b\in \oplus_{i=0}^{r-2}H^{n-m-2i-2}(N,\mz_2)\overline{\xi}^i \subset H^m(P(\eta_N),\mz_2)$.
\end{prop}

%In section 3.3, combing Proposition 1.2, Theorem 1.3 and Wu formula, we will obtain:
%\begin{coro} The Stiefel-Whitney class of $\tilde{M}$ has the form:
%\[
%w(\tilde{M})=\pi^*w(M )+              %\alpha_{j_P}^{*}(\frac{1}{\xi}(p^*w(N)(\sum_{i=0}^{r}w_{2i}(\eta_N)(1+\xi)^{r+1-i})-i_{N}^{*}w(M)) )
%\]
%\end{coro}

In section 4, we want to apply some properties of the Wu classes of topological blow ups to a geometric problem, which is also a partial motivation for the author to discuss the Wu classes of topological blow ups.

Let $F_d\hookrightarrow \mc P^{n+1}$ be a smooth hypersurface of degree $d$, $n$ odd. In his paper \cite{Browder2}, Browder found an obstruction $v_{n+1}(W)\in H^{n+1}(\mc P^{n+1},\mz_2)$, where $W$ is the disk bundle of a certain vector bundle over $\mc P^{n+1}$, such that:\\
(1). If $v_{n+1}(W)\neq 0$, then there exists an embedded sphere $i:S^n \hookrightarrow F_d$ such that $i_* [S^n]=0\in H^n (F_d,\mz_2)$ and the normal bundle of $S^n$ is not trivial but stable trivial.\\
(2). If $v_{n+1}(W)=0$, then a quadratic form $\psi: H_n (F_d, \mz_2)\longrightarrow \mz_2$, which is used to define the Kervaire invariant of $F_d$, is well defined and for any embedded sphere $j:S^n \hookrightarrow F_d$, $\psi (j_*[S^n])=0 \Longleftrightarrow$ the normal bundle of this sphere is trivial.

In our case, if we can blow up $F_d$ along a submanifold $N$ to obtain $\widetilde{F_d}$, then we can also blow up $\mc P^{n+1}$ along $N$ to obtain $\widetilde{\mc P^{n+1}}$ and an embedding:
\[
\widetilde{F_d} \hookrightarrow \widetilde{\mc P^{n+1}}
\]
Likewise, we can also find an obstruction $v_{n+1}(\tilde{W})\in H^{n+1}(\widetilde{\mc P^{n+1}},\mz_2)$. Using the properties of the Steenrod operations and Wu classes of topological blow ups, we will try to describe this obstruction and find out when this obstruction is vanished.

For the detail of this obstruction and more general discussion, see section 4 and the original article of Browder \cite{Browder2}.

\section{Degree one map and topological blow up}

\subsection{Degree one map}
Let $M$ and $N$ be $n$-dimensional, smooth, closed, $R$-oriented manifolds, here $R$ is $\mz$, ${\mz}_2$ or $\mathbb{Q}$. A differential map $
f: M \longrightarrow N$ is called \textbf{a degree one map} if:
\[
f_* [M]=[N]
\]
where $[M]\in H^n (M,R),\ [N]\in H^n (N,R)$ are the fundamental classes.

We introduce the Gysin map, which will be often used in this paper. Given a differentiable map $h:M_1 \longrightarrow M_2$, where $M_1$ and $M_2$ are two smooth, closed, $R$-oriented manifolds with dimension $m_1,\ m_2$. We define the Gysin map $\alpha_{h}^{*}: H^{*}(M_1,R)\rightarrow H^{*+ m_2 - m_1}(M_2,R)$ as follows:
\[
\begin{CD}
\alpha_{h}^{i}:@. H^i (M_1,R)        @>>>     H^{i+m_2 - m_1}(M_2,R)  \\
      @.             @VV{PD}V                   @AA{PD^{-1}}A   \\
               @. H_{m_1 - i}(M_1,R) @>h_{*}>>   H_{m_2 -i} (M_2,R)
\end{CD}
\]
here $PD$ is the Poincar\'{e} duality map. More precisely,
$$\alpha_{h}^{*}:= PD^{-1}h_* PD$$

Next, we discuss the properties of $\alpha_{f}^{*}$ when $f:M^n \longrightarrow N^n$ is a degree one map. In this case, we see the map $\alpha_{f}^{*}: H^{*}(M,R)\rightarrow H^{*}(N,R)$ preserves the degree of the cohomology group. Let $f^*: H^{*}(N,R)\rightarrow H^{*}(M,R)$ be the ring homomorphism induced by $f$, we can compose $f^*$ and $\alpha_{f}^{*}$ to get a homomorphism $H^{*}(N,R)\rightarrow H^{*}(N,R)$. We have:
\begin{lem}
Let $f:M^n \longrightarrow N^n$ be a degree one map, then $\alpha_{f}^{*}f^* = id: H^{*}(N,R)\rightarrow H^{*}(N,R)$
\end{lem}
\begin{proof}
For any $x\in H^*(N,R)$, we have identity:
\[
f_* (PD(f^*(x)))=f_* (f^*(x) \cap [M] )= x\cap f_* [M]=x\cap [N] = PD(x)
\]
according to the definition of $\alpha_{f}^{*}$, $x=PD^{-1}f_* (PD(f^*(x)))=\alpha_{f}^{*}f^*(x)$.
\end{proof}
This identity $\alpha_{f}^{*}f^* = id$ implies the split decomposition of $H^* (M,R)$, namely,
$$H^* (M,R)= f^* H^* (N,R)\oplus \text{Ker} \alpha_{f}^{*}$$
Furthermore, the cup product on $H^*(M.R)$ admits a pairing:
\[
H^p(M,R)\otimes H^{n-p} (M,R) \longrightarrow R
\]
\[
x\otimes y \longmapsto <x\cup y,[M]>
\]
which also induces a pairing:
\[
\text{Ker}\alpha^{p}_{f} \otimes \text{Ker}\alpha^{n-p}_{f} \longrightarrow R
\]

The following proposition describe the behavior of $f^*H^*(N)$ and Ker$\alpha_{f}^{*}$ under these pairings.
\begin{prop}
Under the paring $H^p(M,R)\otimes H^{n-p} (M,R) \longrightarrow R$, $Ker\alpha^{p}_{f}$ is orthogonal to
$f^*(H^{n-q}(N))$. And the pairing $Ker\alpha^{p}_{f} \otimes Ker\alpha^{n-p}_{f} \longrightarrow R$ is nonsingular when $R$ is a field (${\mz}_2, \ \mathbb{Q}$). If $R=\mz$, it is nonsingular on ($Ker\alpha^{p}_{f} /torsion )\otimes (Ker\alpha^{n-p}_{f}/torsion) \longrightarrow \mz$
\end{prop}
\begin{proof}
We refer to \cite{Browder} I.2.9 page 10.
\end{proof}

\subsection{Topological blow up}
Let $M$ be a smooth, closed, $\mz$-oriented manifold with dimension $n$ and $i_N:N \hookrightarrow M$ be a closed submanifold of $M$ with codimension $2r > 2$. If the normal bundle $\eta_N$ of $N$ admits a complex structure, we can blow up $M$ along $N$ as follows:

Step 1: Since the normal bundle $\eta_N$ is a complex vector bundle, we can do the complex projectivization of $\eta_N$ to get $P(\eta_N)$, which is a fiber bundle over $N$ with fiber $\mc P^{r-1}$:
\[
\begin{CD}
\mc P^{r-1} @>>> P(\eta_N)\\
@.                 @VV p V \\
            @.      N
\end{CD}
\]

Step 2: There is a tautological complex line bundle $S$ over $P(\eta_N)$ such that $S$ is the subbundle of $p^* \eta_N$ defined by:
\[
S:= \{(l,v)\in p^* \eta_N \subset P(\eta_N)\times \eta_N |\  v \in l \}
\]
see the following diagram:
\[
\begin{CD}
S         @>>>   p^* \eta_N @>>> \eta_N \\
@VVV             @VVV             @VVV  \\
P(\eta_N) @=     P(\eta_N)  @>p>>  N
\end{CD}
\]
and there is a diffeomorphism $\phi$ between $S-P(\eta_N)$ and $\eta_N - N$:
\[
\phi: S-P(\eta_N) \longrightarrow \eta_N - N
\]
\[
(l,v)\longmapsto v
\]

Step 3: Use this diffeomorphism $\phi$, we construct a closed, smooth, $\mz$-oriented manifold $\tilde{M}$:
\[
\tilde{M}:= (M-N) \sqcup S / x \sim \phi(x), \ x\in S-N
\]

\begin{defn}
$\tilde{M}$ is called the topological blow up of $M$ along $N$.
\end{defn}

\begin{rem}
For any closed submanifold $N\hookrightarrow M$ with normal bundle $\eta_N$, we can do the real projectivization $P_{\mathbb{R}}(\eta_N)$ and we can also construct the real blow up ${\tilde{M}}_{\mathbb{R}}$ of $M$ along $N$.
\end{rem}

\begin{prop}
There exists a map $\pi: \tilde{M} \rightarrow M$ with:\\
(1) $\pi: {\pi}^{-1}(M-N) \longrightarrow M-N$ is a diffeomorphism.\\
(2) $\pi: {\pi}^{-1}(\eta_N)=S \rightarrow \eta_N$ coincides with the composition $S\hookrightarrow p^*\eta_N \rightarrow \eta_N$. In particular, $\pi^{-1}(N)=P(\eta_N)$ is a submanifold of $\tilde{M}$ with normal bundle $S$. \\
(3) $\pi$ is a degree one map for the coefficients $\mz, \ \mz_2,\ \mathbb{Q}$.
\end{prop}
\begin{proof}
(1) and (2) are induced directly from the construction of $\tilde{M}$, we only prove (3). We first consider the case of coefficient $\mz$. Choose a point $x$ and a small ball $D_x$ such that $x\in D_x \subset M-N \cong \tilde{M}- P(\eta_N)$. We see:
\[
\begin{CD}
H_n (\tilde{M},\mz) @>\cong>> H_n (\tilde{M}, \tilde{M}-x,\mz) @= H_n (D_x,D_x -x,\mz)\\
@VV \pi V                   @VVV                                    @|  \\
H_n (M,\mz)         @>\cong>> H_n (M, M-x,\mz)                           @= H_n (D_x,D_x -x,\mz)
\end{CD}
\]
This diagram implies $(\pi)_* ([\tilde{M}])=[M]$. The cases of the coefficients $\mz_2$ and $\mathbb{Q}$  follow immediately.
\end{proof}

We continue to discuss the cohomology group of $H^*(\tilde{M},R)$. The cohomolgy of the blow up has already been researched in \cite{Duan},\cite{Griff},\cite{Macduff}. But in this paper, we want to make it under the point of view of degree one map, which will be useful in our further discussions.

We first discuss the coefficient $\mz$. Since $\pi: \tilde{M} \longrightarrow M$ is a degree one map and we have direct decomposition $H^*(\tilde{M},\mz)=\pi^* H^*(M,\mz) \oplus \text{Ker}\alpha_{\pi}^*$.
In the rest part of this subsection, we'll mainly describe $\text{Ker}\alpha_{\pi}^*$.
We have exact sequences:
\[
\begin{CD}
\cdot\cdot\cdot @>>>H^i(\tilde{M},\tilde{M}-P(\eta_N),\mz) @>r_{\tilde{M}}>>H^i(\tilde{M},\mz) @>>> H^i (\tilde{M}-P(\eta_N),\mz)@>>>\cdot\cdot\cdot\\
 @.                 @AA \pi^* A                        @AA\pi^*A               @|     @.\\
\cdot\cdot\cdot @>>>H^i(M,M-\eta_N,\mz) @>r_{M}>>H^i(M,\mz) @>>> H^i (M-\eta_N,\mz)@>>>\cdot\cdot\cdot\\
\end{CD}
\]
and commutative diagram:
\[
\begin{CD}
H^i (S,S-P(\eta_N),\mz) @>\cong>> H^i(\tilde{M},\tilde{M}-P(\eta_N),\mz)\\
@AA\cup U_S A                                   @Vr_{\tilde{M}}VV  \\
H^{i-2}(P(\eta_N),\mz)  @>\alpha_{j_P}^{*}>>        H^i(\tilde{M},\mz)\\
\end{CD}
\]
where $U_S$ is the Thom class of the normal bundle $S$ over $P(\eta_N)$ and $\alpha_{j_P}^{*}$ is the Gysin map induced by the embedding $j_P: P(\eta_N)\hookrightarrow \tilde{M}$.

On the other hand, by Leray-Hirsh theorem, we have isomorphism:
\[
           H^*(N,\mz)[t]/(\sum_{i=0}^{r} c_{r-i}(\eta_N)t^i) \rightarrow H^* (P(\eta_N),\mz)
\]
\[
               \sum m_i t^i \longmapsto \sum (p^*m_i ){\xi}^i
\]
where $-\xi$ is the Euler class of the complex line bundle $S$ and $c_i(\eta_N)$ is the i-th Chern class of the complex bundle $\eta_N$. From now on, for convenience, we identify $H^*(P(\eta_N),\mz)$ with $\oplus_{i=0}^{r-1}H^*(N,\mz)\xi^{i}$.

The main result in this section is:

\begin{prop}
$Ker\alpha_{\pi}^{*}= \alpha_{j_P}^{*}(\oplus_{i=0}^{r-2}H^*(N,\mz){\xi}^i)$
\end{prop}
\begin{proof}
First, we have commutative diagram:
\[
\begin{CD}
H^i(\tilde{M},\mz) @>j_{P}^{*}>>  H^{i}(P(\eta_N),\mz)\\
@AA\alpha_{J_P}^*A                 @| \\
H^{i-2}(P(\eta_N),\mz) @>\cup \xi >>    H^i(P(\eta_N),\mz)
\end{CD}
\]
We see $\cup \xi: \oplus_{i=0}^{r-2}H^*(N,\mz){\xi}^i \longrightarrow  H^*(P(\eta_N),\mz) $ is injective and this implies $\alpha_{J_P}^*$ is injective when it restricts to $\oplus_{i=0}^{r-2}H^*(N,\mz){\xi}^i$.

Second, we need to prove $\alpha_{j_P}^{*}(\oplus_{i=0}^{r-2}H^*(N,\mz){\xi}^i)\subset$ Ker$\alpha_{\pi}^{*}$. %which is equivalent to prove
%$$\alpha_{\pi}^{*}\alpha_{j_P}^{*}(\oplus_{i=0}^{r-2}H^*(N){\xi}^i)=0$$
For any $(p^*a)\xi^{i}\in H^*(P(\eta_N),\mz):$
\[
p_* (PD((p^*a)\xi^i))=p_*((p^*a)\xi^{i}\cap [P(\eta_N)])=p_* (p^* a\cap (\xi^i \cap [P(\eta_N)]))=a\cap p_*(\xi^i\cap [P(\eta_N)])
\]
Since deg($\xi^i\cap [P(\eta_N)])=n-2r+2r-2-2i=2(r-1-i)+$dim$N$, we see $p_*(\xi^i\cap [P(\eta_N)])=0$, when $i<r-1$. This implies $\alpha_{p}^{*}(\oplus_{i=0}^{r-2}H^*(N,\mz){\xi}^i)=PD^{-1}p_*PD(\oplus_{i=0}^{r-2}H^*(N,\mz){\xi}^i)=0$.
Since $\pi \circ j_P =i_N \circ p$,
\[
\begin{CD}
P(\eta_N) @>j_P>> \tilde{M}\\
@VpVV              @V\pi VV\\
N @>i_N >>              M
\end{CD}
\]
we have $\alpha_{\pi}^{*}\alpha_{j_P}^{*}(\oplus_{i=0}^{r-2}H^*(N,\mz){\xi}^i)=\alpha_{i_N}^{*}
\alpha_{p}^{*}(\oplus_{i=0}^{r-2}H^*(N,\mz){\xi}^i)=0$.

Third, we want to prove $\alpha_{j_P}^{*}(\oplus_{i=0}^{r-2}H^*(N,\mz){\xi}^i) + \pi^* H^*(M,\mz)=
\text{Ker}\alpha_{\pi}^{*}\oplus \pi^* H^*(M,\mz) = H^*(\tilde{M},\mz)$, which implies
\[
\alpha_{j_P}^{*}(\oplus_{i=0}^{r-2}H^*(N,\mz){\xi}^i)= \text{Ker}\alpha_{\pi}^{*}.
\]
By chasing the diagram below:
\[
\begin{CD}
\cdot\cdot\cdot @>>>H^i(\tilde{M},\tilde{M}-P(\eta_N),\mz) @>r_{\tilde{M}}>>H^i(\tilde{M},\mz) @>>> H^i (\tilde{M}-P(\eta_N),\mz)@>>>\cdot\cdot\cdot\\
 @.                 @AA \pi^* A                        @AA\pi^*A               @|     @.\\
\cdot\cdot\cdot @>>>H^i(M,M-\eta_N,\mz) @>r_M>>H^i(M,\mz) @>>> H^i (M-\eta_N,\mz)@>>>\cdot\cdot\cdot\\
\end{CD}
\]
we see $H^*(\tilde{M},\mz)=\pi^*{H^*(M,\mz)}+r_{\tilde{M}}(H^*(\tilde{M},\tilde{M}-P(\eta_N),\mz))$. We also know
$r_{\tilde{M}}(H^*(\tilde{M},\tilde{M}-P(\eta_N),\mz))=\alpha_{j_P}^*(H^*(P(\eta_N),\mz))$.
From the commutative diagram:
\[
\begin{CD}
H^{i-2}(P(\eta_N),\mz) @>\cong >>  H^i (S,S-P(\eta_N),\mz) @>\cong >> H^i (\tilde{M}, \tilde{M}-P(\eta_N),\mz)\\
@A\cup e((p^*\eta_N)/S)AA      @AAA                                @A \pi^* AA\\
H^{i-2r}(N,\mz)   @>\cong >> H^i (\eta_N, \eta_N-N ,\mz)  @>\cong>>   H^i (M,M-N,\mz)
\end{CD}
\]
here $e((p^*\eta_N)/S)=\sum_{i=0}^{r-1}c_i(\eta_n)\xi^i=\xi^{r-1}+\sum_{i<r-1} c_i(\eta_N) \xi^i$ is the Euler class of the quotient bundle $p^*\eta_N/S$ over $P(\eta_N)$, we get decomposition
$$H^*(P(\eta_N),\mz)=H^*(N,\mz)e((p^*\eta_N)/S) \oplus (\oplus_{i=0}^{r-2}H^*(N,\mz){\xi}^i)$$
Then, we have:
\[
\alpha_{j_P}^*(H^*(P(\eta_N),\mz))
=\alpha_{j_P}^{*}(\oplus_{i=0}^{r-2}H^*(N,\mz){\xi}^i) + \alpha_{j_P}^*(H^*(N,\mz)e((p^*\eta_N)/S))
\]
\[
\alpha_{j_P}^*(H^*(N,\mz)e((p^*\eta_N)/S))=\pi^*\alpha_{i_N}^{*}H^*(N,\mz)\subset \pi^* H^*(M,\mz)
\]
Thus, we obtain:
\[
H^*(\tilde{M},\mz)=\pi^*{H^*(M,\mz)} + \alpha_{j_P}^{*}(\oplus_{i=0}^{r-2}H^*(N,\mz){\xi}^i)
\]
\end{proof}
The discussion of $H^*(\tilde{M},\mz_2)$ and $H^*(\tilde{M},\mathbb{Q})$ are similar to $H^*(\tilde{M},\mz)$. We have:
\begin{coro}
$$H^*(\tilde{M},\mz_2)=\pi^* H^*(M,\mz_2)\oplus \alpha_{j_P}^{*}(\oplus_{i=0}^{r-2}H^*(N,\mz_2){\overline{\xi}}^i)$$
and
$$H^*(\tilde{M},\mathbb{Q})=\pi^*H^*(M,\mathbb{Q})\oplus \alpha_{j_P}^{*}(\oplus_{i=0}^{r-2}H^*(N,\mathbb{Q})\xi_{Q}^i)$$
where $\overline{\xi}$ and $\xi_{Q}$ is the image of $\xi$ in $H^2(P(\eta_N),\mz_2)$ and $H^2(P(\eta_N),\mathbb{Q})$.
\end{coro}

\subsection{L-class of the topological blow up}
As a quick application, in this subsection, we give a formula of the L-class $L_i (\tilde{M})$ of the topological blow up of $M^n$ along its submanifold $N$ with codimension $2r$ when $4i < min\{2r,\frac{n-1}{2}\}$. In this subsection, the coefficient is $\mathbb{Q}$. Following the notation of subsection 2.1 and 2.2, our result is:
\begin{prop}
Suppose $\tilde{M}$ is the topological blow up of $M$ along the submanifold $N$, when $4i < min\{2r,\frac{n-1}{2}\}$,
$$L_i(\tilde{M})=\pi^{*}L_i(M)+ \alpha_{j_P}^{*}(\xi_{Q}^{-1}(L_i(S)-i_{N}^{*}(L_i (M))))$$
\end{prop}
\begin{proof}
According to Corollary 2.7.
$$L_i(\tilde{M})=\pi^{*}a+ \alpha_{j_P}^{*}(b)\in \pi^*{H^{4i}(M)} \oplus \alpha_{j_P}^{*}(\oplus_{i=0}^{r-2}H^*(N)\xi_{Q}^i)  $$
When $4i < min\{2r,\frac{n-1}{2}\}$, for any smooth map $f: M \longrightarrow S^{n-4i}$, since dim$N=n-2r<n-4i$, we can choose a regular point $y\in S^{n-4i}-f(N)$ to get the smooth manifold $f^{-1}(y)\subset M-N=\tilde{M}-P(\eta_N)$.
\[
\begin{CD}
\tilde{f}^{-1}(y) @>>> \tilde{M}-P(\eta_N)\subset \tilde{M} @>\tilde{f}=\pi f>> S^{n-4i} \\
 @|                                 @VV\pi V                              @|\\
 f^{-1}(y)  @>>>    M-N\subset M  @>f>> S^{n-4i}
\end{CD}
\]
Let $u\in H^{n-4i}(S^{n-4i})$ be the generator. By the property of $L$ class (for example \cite{Milnor}), we have
\begin{equation*}\label{xx}
\begin{split}
<L_{i}(\tilde{M})\cup \tilde{f}^{*}u,[\tilde{M}]> & = \sigma (\tilde{f}^{-1}(y))\\
                                                     & = \sigma (f^{-1}(y))\\
                                        & = <L_i (M)\cup f^{*}u,[M]>\\
\end{split}
\end{equation*}
where $\sigma()$ denotes the index of a closed manifold. On the other hand, by proposition 2.2, ker$\alpha_{\pi}^*$ is orthogonal to $\pi^* H^*(M)$ under the pairing of cup product, we get:
\begin{equation*}\label{xx}
\begin{split}
<L_{i}(\tilde{M})\cup \tilde{f}^{*}u,[\tilde{M}]> & = <\pi^{*}a \cup \tilde{f}^{*}u,[\tilde{M}]>\\
                                                     & = <\pi^*(a\cup f^{*}u),[\tilde{M}]>\\
                                        & = <a\cup f^{*}u,[M]>\\
\end{split}
\end{equation*}
We know that the equation $<L_i (M)\cup f^{*}u,[M]>=\sigma (f^{-1}(y))$ for every map $f:M \longrightarrow S^{n-4i} $ determines $L_i(M)$ when $i< \frac{n-1}{2}$ (\cite{Milnor}), then we get $a=L_i (M)$ and
\[
 L_i (\tilde{M})= \pi^{*}L_i(M) + \alpha_{j_P}^{*}(b)
\]
For the normal bundle $ S \subset \tilde{M}$, we have
\begin{equation*}\label{xx}
\begin{split}
L_i(S)                  &= L_i (\tilde{M})|S       \\
                        &=j_{P}^{*}L_i (\tilde{M}) \\
                        & = j_{P}^{*}\pi^{*}L_i(M) + b\xi_Q \\
                        & = i_{N}^{*}L_i(M)+ b \xi_Q \\
\end{split}
\end{equation*}
Since $b\xi_Q \in \sum_{i=1}^{r-1}H^*(N)\xi_{Q}^{i}$, we have:
\[
b=\xi_{Q}^{-1}(L_i(S)-i_{N}^{*}L_i(M))
\]
Furthermore, the tangent bundle of $S$ is a complex vector bundle with Chern class:
\[
c(S)=(1-\xi_Q)(\sum_{i=0}^{r}c_i(\eta_N)(1+\xi_Q)^{r-i}) \in H^*(S,\mz)=H^*(P(\eta_N),\mz)
\]
$L_i(S)$ can be calculated by the L-polynomial $L(p_1(S),p_2(S),\cdot\cdot\cdot)$ with $ \sum_i (-1)^i p_i(S)=c(S)c(\overline{S})$.

\end{proof}

\section{Steenrod operations and Wu classes of the topological blow up}
In this section, all the coefficients of the cohomology groups are ${\mz}_2$.
\subsection{Steenrod operations on the topological blow up}
Let $\tilde{M}$ be the topological blow up of $M$ along the submanifold $N$, we first determine the Steenrod operations on $H^* (\tilde{M})$. Our result is:
\begin{prop}
For any element $x=\pi^* (x_1)+ \alpha_{j_P}^{*}(x_2) \in \pi^*H^*(M)\oplus \alpha_{j_P}^{*}(\oplus_{i=0}^{r-2}H^*(N){\overline{\xi}}^i)=H^*(\tilde{M})$:
\[
Sq(x)=\pi^* Sq(x_1) + \alpha_{j_P}^{*}((1+\overline{\xi})Sq(x_2))
\]
\end{prop}
\begin{proof}
We only have to prove $Sq(\alpha_{j_P}^{*}(x_2))=\alpha_{j_P}^{*}((1+\overline{\xi})Sq(x_2))$. Since we have commutative diagram:
\[
\begin{CD}
H^* (S,S- P(\eta_N)) @>\cong >>   H^*(\tilde{M},\tilde{M}-P(\eta_N))\\
@AA\cup \overline{U_S}A                             @V r_{\tilde{M}}VV\\
H^*(P(\eta_N))       @>\alpha_{j_P}^{*}>>  H^* (\tilde{M})
\end{CD}
\]
where $\overline{U_S}$ is the mod 2 Thom class of the normal bundle $S$ of $P(\eta_N)$ in $\tilde{M}$. The Stiefel-Whitney class $w(S)$ of $S$ is equal to $1+\overline{\xi}$ and we have:
\begin{equation*}\label{xx}
\begin{split}
Sq(\alpha_{j_P}^{*}(x_2))&= Sq( r_{\tilde{M}}(x_2 \cup \overline{U_S}))\\
                         &= r_{\tilde{M}}(Sq(x_2\cup \overline{U_S}))\\
                         &= r_{\tilde{M}}(Sq(x_2)\cup Sq(\overline{U_S}))\\
                         &= r_{\tilde{M}}(Sq(x_2)\cup w(S)\cup \overline{U_S})\\
                         &= r_{\tilde{M}}((Sq(x_2)\cup(1+\overline{\xi})) \cup \overline{U_S})\\
                         &= \alpha_{j_P}^{*}((1+\overline{\xi})Sq(x_2))\\
\end{split}
\end{equation*}
here the identity $Sq(\overline{U_S})=w(S)\cup \overline{U_S}$ is just Thom's formula of the Stiefel-Whitney characteristic classes $w(S)=1+\overline{\xi}$.
\end{proof}

\begin{coro}
For any element $x=\pi^* (x_1)+ \alpha_{j_P}^{*}(x_2) $ and ${Sq}^i$, we have:
$${Sq}^i(x)=\pi^* {Sq}^i (x_1) + \alpha_{j_P}^{*}(Sq^i x_2 + \overline{\xi} Sq^{i-2}(x_2))$$
\end{coro}
\begin{rem}
The element $\alpha_{j_P}^{*}((1+\overline{\xi})Sq(x_2))$, it is not necessarily belong to Ker$ \alpha_{j_P}^{*}$.
\end{rem}

\subsection{Wu classes of the blow up}
Given a smooth manifold $(M^n,\partial M)$, the Steenrod operator $Sq=\sum_{i=0} {Sq}^i:H^*(M,\partial M)\rightarrow H^*(M,\partial M)$ determines a linear form on $H^*(M,\partial M)$:
\[
H^*(M,\partial M )\longrightarrow {\mz}_2
\]
\[
x \mapsto <Sq(x),[M]>
\]
where $[M]\in H_n(M,\partial M)$ is the fundamental class of the Poincar\'{e} pair $(M,\partial
M)$. Since the cup product induces the isomorphism $H^*(M)\cong Hom(H^*(M,\partial M),{\mz}_2)$, there exists a unique element $v(M)=1+v_1(M)+v_2(M)+\cdot\cdot\cdot\in H^*(M) $ such that for each $x\in H^*(M,\partial M)$:
\[
<v(M)\cup x,[M]>= <Sq(x),[M]>
\]
\begin{defn}
$v(M)=\sum_{i=0}v_i (M)$ is called the Wu class of $M$.
\end{defn}
Since Wu formula will be often used, we denote
$$\chi(Sq): = {Sq}^{-1}: H^*(M)\longrightarrow H^*(M)$$
by the inverse of the isomorphism $Sq:H^*(M)\rightarrow H^*(M)$, we recall:
\begin{thm}[Wu formula]
$(M,\partial M)$ is a smooth compact manifold, then we have
\[
Sq(v(M))= w(M)
\]
or equivalently
\[
v(M)=\chi(Sq)(w(M))
\]
where $w(M)$ is the Stiefel-Whitney class of $M$.
\end{thm}

Wu class keeps the decomposition under the degree one map. Indeed, let $M_1$ and $M_2$ be two smooth manifolds and $f:M_1\longrightarrow M_2$ be a degree one map, we have
\begin{lem}
$v(M_1)=f^* v(M_2) + \overline{v}\in f^*H^*(M_2)\ \oplus$ $Ker\alpha_{f}^{*}=H^*(M_1)$.
\end{lem}
\begin{proof}
Since $f^*H^*(M_2)\ \oplus$ Ker$\alpha_{f}^{*}=H^*(M_1)$, suppose $v(M_1)=f^* a + \overline{v}$. We want to prove $a=v(M_2)$. For every $x\in H^*(M_2)$, we see
\begin{equation*}
\begin{split}
<v(M_1)\cup f^* x,[M_1]> & = <f^*a\cup f^* x,[M_1]>\\
                                     & = <f^*(a\cup x),[M_1]> \\
                                     & = <a\cup x, [M_2]>
\end{split}
\end{equation*}
On the other hand:
\begin{equation*}
\begin{split}
<v(M_1)\cup f^* x,[M_1]> & = <Sq(f^*x),[M_1]>\\
                                     & = <f^*(Sq(x)),[M_1]> \\
                                     & = <Sq(x), [M_2]>\\
                                     & = <v(M)\cup x,[M_2]>
\end{split}
\end{equation*}
Then we get $<a\cup x,[M_2]>=<v(M_2)\cup x,[M_2]>$ for every $x\in H^*(M_2)$. Since the coefficients of the cohomology group is ${\mz}_2$, we conclude: $a=v(M_2)$.
\end{proof}

Next we discuss the Wu classes of the topological blow up. Suppose $\tilde{M}$ is the topological blow up of $M$ along the submanifold $N$. There is a natural degree one map $\pi:\tilde{M}\longrightarrow M$. According to Lemma 3.6. and corollary 2.7, we have
\begin{coro}
$v(\tilde{M})=\pi^*(v(M)) + \alpha_{j_P}^{*}(b)$, where $b\in \sum_{i=0}^{r-2}H^*(N)\overline{\xi}^i$.
\end{coro}

We know $S$ is the normal bundle of $P(\eta_N)$ and the closure $\overline{S}$ is a submanifold with boundary of $\tilde{M}$. Denote $v(\overline{S})\in H^*(\overline{S})=H^*(P(\eta_N))$ by the Wu class of $(\overline{S},\partial \overline{S})$, we have:
\begin{prop}
$b={\overline{\xi}}^{-1}(v(\overline{S})-{i_N}^* v(M))$
\end{prop}
\begin{proof}
Since $T\tilde{M}{|}_{\overline{S}}=T\overline{S}$, we have $w(\overline{S})=j_{P}^{*}w(M)\in H^*(P(\eta_N))$. By the
Wu formula:
\begin{equation*}
\begin{split}
v(\overline{S}) &= \chi(Sq)w(\overline{S})\\
                &= \chi(Sq)j_{P}^{*}w(\tilde{M}) \\
                & = j_{P}^{*} \chi(Sq)(w(\tilde{M}))\\
                & = j_{P}^{*} v(\tilde{M})
\end{split}
\end{equation*}
On the other hand,
\[
j_{P}^{*} (v(\tilde{M}))=j_{P}^{*}(\pi^*v(M) + \alpha_{j_P}^{*}(b))=i_{N}^{*}v(M) + b\overline{\xi}
\]
we get
\[
b\overline{\xi} = v(\overline{S})-{i_N}^* v(M) \in \sum_{i=1}^{r-1}H^*(N)\overline{\xi}^i
\]
Since $\cup\overline{\xi}: \sum_{i=0}^{r-2}H^*(N)\overline{\xi}^i\longrightarrow \sum_{i=1}^{r-1}H^*(N)\overline{\xi}^i$ is bijective, its inverse ${\overline{\xi}}^{-1}:\sum_{i=1}^{r-1}H^*(N)\overline{\xi}^i\longrightarrow \sum_{i=0}^{r-2}H^*(N)\overline{\xi}^i$ is well-defined, we get:
\[
b= {\overline{\xi}}^{-1}(v(\overline{S})-{i_N}^* v(M))
\]
\end{proof}

We continue to describe $v(\overline{S})$, our result is:
\begin{prop}
$v(\overline{S})=v(P(\eta_N))\chi(Sq)(1+\overline{\xi})\in H^*(P(\eta))$
\end{prop}
\begin{proof}
Since $S$ is the normal bundle of $P(\eta_N)$ and let $\pi_{P}$ be the projection map. We see $TS=\pi_{P}^{*}TP(\eta_N)\oplus \pi_{P}^{*}S$ and
$w(\overline{S})=(1+\overline{\xi})w(P(\eta_N))$
By Wu formula, we get: $v(\overline{S})= v(P(\eta_N))\chi(Sq)(1+\overline{\xi})$.

%By the definition of Wu class, we need to detect the homomorphism $Sq: H^*(\overline{S},\partial \overline{S}) %\longrightarrow H^*(\overline{S},\partial \overline{S}).$ Observe that
%\[
%\cup U_S: H^*(P(\eta_N)) \cong
%H^*(\overline{S},\overline{S}-P(\eta_N))\cong H^*(\overline{S},\partial \overline{S})
%\]
%For every element $x\cup U_2 \in H^*(\overline{S},\partial\overline{S})$, we have:
%\[
%Sq(x\cup U_S)=Sq(x)\cup Sq(U_S)=(1+\xi)Sq(x)\cup U_S
%\]
%Thus, we get a diagram:
%\[
%\begin{CD}
%H^*(\overline{S},\partial \overline{S}) @>Sq>> H^*(\overline{S},\partial \overline{S})\\
% @AA\cup U_SA                           @AA\cup U_S A \\
% H^*(P(\eta_N))  @>(1+\xi)Sq >>         H^*(P(\eta_N))
%\end{CD}
%\]
%Furthermore, for every element $x\cup U_S \in H^*(\overline{S},\partial\overline{S})$
%\begin{equation*}
%\begin{split}
%<Sq(x\cup U_S),[\overline{S}]>   & = <(1+\xi)Sq(x)\cup U_S,[\overline{S}]>\\
%                                 & = <(1+\xi)Sq(x),[P(\eta_N)]>\\
%                                 & = <Sq(\chi(Sq)(1+\xi)x),[P(\eta_N)]>\\
%                                 & = <v(P(\eta_N))\chi(Sq)(1+\xi)x,[P(\eta_N)]> \\
%                                 & = <v(P(\eta_N))\chi(Sq)(1+\xi)(x\cup U_S), [\overline{S}]>
%\end{split}
%\end{equation*}
%By the definition of Wu class, we obtain
\end{proof}

\begin{lem}
$v(P(\eta_N))=p^* (v(N)) \chi(Sq)(\sum_{i=1}^{r}w_{2i}(\eta_N)(1+\overline{\xi})^{n+1-i})$
\end{lem}
\begin{proof}
Since $P(\eta_N)$ is the complex projectivization of $\eta_N$, by the formula $TP(\eta_N)=p^* TN \oplus Hom(S,\pi^*\eta_N /S)$, we have:
\begin{equation*}
\begin{split}
w(P(\eta_N)) & =p^*(w(N)) w(Hom(S,\pi^*\eta_N /S))\\
             & =p^*(w(N))w(p^*\eta_N \otimes S^{-1})\\
             & =p^*(w(N))(\sum_{i=0}^{r}w_{2i}(\eta_N)(1+\overline{\xi})^{r-i})
\end{split}
\end{equation*}
here $S^{-1}=Hom(S,\mc)$ is the inverse of the complex line bundle $S$.
Compose $\chi(Sq)$ to both sides and use the Wu furmula:
\begin{equation*}
\begin{split}
v(P(\eta_N)) &=\chi(Sq)(w(P(\eta_N)))\\
             & =p^*(\chi(Sq)w(N))\chi(Sq)(\sum_{i=0}^{r}w_{2i}(\eta_N)(1+\overline{\xi})^{r-i})\\
             &= p^*v(N)\chi(Sq)(\sum_{i=0}^{r}w_{2i}(\eta_N)(1+\overline{\xi})^{r-i})\\
\end{split}
\end{equation*}
\end{proof}

Finally, combining these lemmas and properties, we get our main theorem:
\begin{thm} The Wu class of the topological blow up $\tilde{M}$ has the form:
\[
v(\tilde{M})=\pi^* v(M)+ \alpha_{j_P}^{*}({\xi}^{-1}(p^*v(N)\chi(Sq)(\sum_{i=0}^{r}w_{2i}(\eta_N)(1+\overline{\xi})^{r+1-i})-i_{N}^{*}(v(M))) )
\]
\end{thm}

\subsection{On the vanishing of Wu class}
We have obtained the formula of Wu class of topological blow ups, however, it seems not easy to calculate. In some cases, we only want to determine whether $v_{m}(\tilde{M})$ is zero or not rather than to determine the element itself.

In this subsection, we will discuss the vanishing condition of Wu class in two cases. One case is the Wu classes $v_m(\tilde{M})$ of the topological blow up $\tilde{M}$. The other case is the Wu classes of a disk bundle over a compact manifold, which will be used in section 4.

Let $\tilde{M}$ be the topological blow up of $M$ along the submanifold $N$.
\begin{prop}
$v_{m}(\tilde{M})=0$ if and only if\\
(1). $v_m (M)= 0$\\
(2). ${Sq}^m b + (Sq^{m-2}b) \overline{\xi}=0$ for all $b\in \oplus_{i=0}^{r-2}H^{n-m-2i-2}(N)\overline{\xi}^i \subset H^{n-m-2}(P(\eta_N))$.
\end{prop}
\begin{proof}
First we suppose $v_m(\tilde{M})=0$. By the definition of $v_m(\tilde{M})$, we see $v_m(\tilde{M})=0$ if and only if the map
\[
Sq^m: H^{n-m}(\tilde{M}) \longrightarrow H^{n}(\tilde{M})={\mz}_2
\]
is zero. Since we have decomposition:
\[
H^{n-m}(\tilde{M})= \pi^* H^{n-m}(M) \oplus \alpha_{j_P}^{*}(\oplus_{i=0}^{r-2}H^{n-m-2i-2}(N)\overline{\xi}^i)
\]
For any $a\in H^{n-m}(M)$, $Sq^m (\pi^* a)=\pi^* (Sq^m a)=0$. Since $\pi^*$ is injective, we have $Sq^m a = 0$, for all $a\in H^{n-m}(M)$. This implies $v_m (M)=0$. For any $b\in \oplus_{i=0}^{r-2}H^{n-m-2i-2}(N)\overline{\xi}^i$, $Sq^m\alpha_{j_P}^{*}(b)=0$. By Corollary 3.2,
\[
Sq^m\alpha_{j_P}^{*}(b) = \alpha_{j_P}^{*}(Sq^m b + (Sq^{m-2}b)\overline{\xi})=0
\]
Since $\alpha_{j_P}^*: H^{n-2}(P(\eta_N))\longrightarrow H^n(M)$ is an isomorphism, we obtain
\[
Sq^m b + (Sq^{m-2}b)\overline{\xi}=0, \ b\in \oplus_{i=0}^{r-2}H^{n-m-2i-2}(N)\overline{\xi}^i
\]

For the "if" part, for any $\pi^*a + \alpha_{j_P}^{*}(b)\in H^{n-m}(\tilde{M})$,
\[
Sq^m(\pi^*a + \alpha_{j_P}^{*}(b)) = \pi^*(Sq^m a)+ \alpha_{j_P}^{*}(Sq^m b + (Sq^{m-2}b)\overline{\xi})=0
\]
This implies $v_m (\tilde{M})=0$.
\end{proof}

\begin{exmp}
We blow up $\mc P^n$ along a point $p$ to obtain $\widetilde{\mc P^n}$. First, $H^{odd}(\mc P^n)=H^{odd}(\tilde{\mc P^n})=0$, then $v_i(\widetilde{\mc P^n})=0, \ i$ odd. Second, we know $H^*(\mc P^n)=\mz_2 [x]/(x^{n+1})$, deg$x=2,$ and $P(\eta_p)\cong\mc P^{n-1}$ with $H^*(P(\eta_p))=\mz_2 [\xi]/(\xi^{n})$,  deg$\xi=2.$

For the generator $x^{n-m}\in H^{2n-2m}(\mc P^n)$,
\[
Sq^{2m} x^{n-m}=
\left( \begin{array}{c}
n-m\\
m
\end{array}
\right)x^n
\]
For the generator $\xi^{n-m-1}\in H^{2n-2m-2}(P(\eta_p))$
\begin{equation*}
\begin{split}
Sq^{2m} \xi^{n-m-1}+ \xi Sq^{2m-2}\xi^{n-m-1} &=  \left(\begin{array}{c}
n-m-1\\
m
\end{array}
\right)\xi^{n-1}+\left(
\begin{array}{c}
n-m\\
m-1
\end{array}
\right)\xi^{n-1}\\
            &= {\left(\begin{array}{c}
n-m\\
m
\end{array}
\right)
\xi^{n-1}}
\end{split}
\end{equation*}

So we get $v_{2m}(\widetilde{\mc P^n}) = 0 $ if and only if $\left(\begin{array}{c}
n-m\\
m
\end{array}
\right)\equiv 0 \ (mod \ 2)$.
\end{exmp}

Next, we will discuss the vanishing of the Wu class of a disc bundle, which will be used in next section.

Let $M$ be a smooth, closed, n-manifold and $E$ be a m-dimensional vector bundle over $M$ with projection map $p$. By choosing a metric $<,>$ on $E$, we denote the disk bundle of $E$ by $D(E):=\{v\in E|<v,v>\leqslant 1\}$ and the sphere bundle of $E$ by $S(E):=\{v\in E|<v,v>=1\}$. $D(E)$ is a compact manifold with boundary $S(E)$. Its Wu class is:
\begin{prop}
$v(D(E))=v(M)\chi(Sq)(w(E))$, $w(E)$ is the Stiefel-Whitney class of $E$.
\end{prop}
\begin{proof}
Since $TD(E)=p^* TM \oplus p^* E$, we get $w(D(E))=w(M)w(E)\in H^*(D(E))=H^*(M)$. By Wu formula,
$v(D(E))=\chi(Sq)(w(D(E)))=\chi(Sq)(w(M))\chi(Sq)(w(E))=v(M)\chi(Sq)(w(E))$.
\end{proof}

We want to know when $v_i(D(E))=0$.
\begin{prop}
$v_i(D(E))=0$ if and only if for any $x\in H^{n-i}(M)$,
\[
<w(E)Sqx,[M]>=0
\]
\end{prop}
\begin{proof}
We want to detect the map $Sq^i: H^{n+m-i}(D(E),S(E))\longrightarrow H^{n+m}(D(E),S(E))=\mz_2$. By Thom isomorphism theorem, we get a commutative diagram:
\[
\begin{CD}
H^{n+m-i}(D(E),S(E)) @>Sq^i>> H^{n+m}(M)\\
 @AA\cup U_E A                           @AA\cup U_E A \\
 H^{n-i}(M)  @>\Phi>>         H^n(M)
\end{CD}
\]
here $U_E$ is the mod 2 Thom class and we know by Thom's formula: $Sq^j U_E= w_j(E)U_E$. For any $x\in H^{n-i}(M)$,
\[
Sq^i(xU_E)=\sum_{j=0}^{i}Sq^j(x)Sq^{i-j}U_E =(\sum_{j=0}^{i}w_{i-j}(E)Sq^j x)U_E
\]
then the map $\Phi$ has the expression $\Phi(x)=\sum_{j=0}^{i}w_{i-j}(E)Sq^j x$. By the commutative diagram, $Sq^i=0 $ if and only if $\Phi=0$. Also $\Phi=0$ is equivalent to:
\[
<\sum_{j=1}^{i}w_{i-j}(E)Sq^j x,[M]>=<w(E)Sqx,[M]>=0
\]
for any $x\in H^{n-i}(M)$.
\end{proof}

\section{Geometric applications}

\subsection{Wu class as a geometric obstruction}

In his paper \cite{Browder2} (1.5) pp 95, Browder proved:
\begin{thm}[Browder] Suppose $M^{2n}\times \mathbb{R}^q \subset W , \ n\neq 1,3$ or $7$, $W$ is 1-connected. $(W,M)$ is n-connected and suppose $v_{n+1}(W)\neq 0$. Then there exists an embedded $S^n \subset M^{2n}$ and $U^{n+1}\subset M^{2n}\times \mathbb{R}^{q+1}$ with $\partial U=S^n$ such that the normal bundle $\xi$ to $S^n$ in $M^{2n}$ is nontrivial, but $\xi+\epsilon^1$ is trivial, where $\epsilon^1$ is the trivial one dimensional real vector bundle. Hence $S^n$ is homologically trivial $(mod \ 2)$ with nontrivial normal bundle.
\end{thm}

On the other hand, if $v_{n+1}(W)=0$, Browder also proved:
\begin{thm}[Browder] Suppose $M^{2n}\times \mathbb{R}^q \subset W$, $n$ odd, $\ n\neq 1,3$ or $7$, $W$ is 1-connected. $(W,M)$ is n-connected and suppose $v_{n+1}(W)=0$. \\
(1). A quadratic form $\psi: K\longrightarrow \mz_2$, where $K=Ker(H_n(M,\mz_2)\rightarrow H_n(W,\mz_2))$, is well defined.\\
(2). For the embedding $\varphi:S^n\hookrightarrow M^{2n}$ with $\varphi$ nullhomotopic in $W$, $\psi(\varphi_*[S^n])=0$ if and only if the normal bundle of $\varphi(S^n)$ is trivial.
\end{thm}
\begin{proof}
We refer to \cite{Browder2} (1.4) pp 93 and (1.7) pp 97.
\end{proof}

For a smooth $F_d$ in $\mc P^{n+1}$, although $(\mc P^{n+1},F_d)$ is n-connected, the normal bundle of $F_d$ is not trivial and we can not apply Browder's theorems to get the obstruction directly. We need a little change. Technically, we need to "thicken" $\mc P^{n+1}$ to make the normal bundle of $F_d$ trivial and use Browder's theorems to find that obstruction.

We consider a more general case. Let $M^{2n+2}$ be a 1-connected, closed, smooth manifold and $i_N :N^{2n}\hookrightarrow M^{2n+2}$ be a closed submanifold of $M^{2n+2}$. Suppose $\pi_i (M,N)=0,\ i<n+1.$

First, we want to "thicken" $M$. Since $H^2 (M,\mz)=[M, BSO(2)]=[M,BU(1)]$, there exists a complex line bundle $L$ over $M$ with the Euler class $e(L)=PD^{-1}((i_N)_*[N])\in H^2 (M,\mz)$. When we pull back $L$ to $N$, the complex line bundle ${i_N}^* L$ is isomorphic to the normal bundle of $N$ because these two real oriennted 2-dimensional vector bundles have the same Euler classes. So we can identify $i_{N}^* L$ to the normal bundle of $N$ in $M$.

Let $-L$ be the vector bundle over $M$ which is stable inverse to $L$, i.e., $L \oplus -L$ is a trivial bundle and let $W=D(-L)$, the disk bundle of $-L$. Then we see the normal bundle of the embedding:
\[
N\hookrightarrow M\hookrightarrow W
\]
is trivial.

Second, by Browder's theorem, we find the obstruction is just $v_{n+1}(W)$. The point is to determine whether this obstruction is zero. More precisely, we have:
\begin{prop}
$v_{n+1}(W)=0$ $\Longleftrightarrow$ for all $x\in H^{n+1}(M^{2n+2},\mz_2)$
\[
<\frac{Sqx}{1+u_N},[M^{2n+2}]> = 0
\]
where $u_N =[e(L)]\in H^2(M,\mz_2)$.
\end{prop}
\begin{proof}
Since $W=D(-L)$, by proposition 3.15, $v_{n+1}(W)=0 \Longleftrightarrow \ <w(-L)Sqx,[M]>=0$ for all $x\in H^{n+1}(M,\mz_2)$. We know $w(L)=1+u_N$ and $w(-L)=\frac{1}{1+u_N}$. Then we have:
\[
<w(-L)Sqx,[M]>=<\frac{Sqx}{1+u_N},[M]>
\]
\end{proof}

\begin{coro}
If $\exists \ x \in H^{n+1}(M,\mz_2)$ such that $<\frac{Sqx}{1+u_N},[M]>\neq 0$, then there exists an embedding $S^n \hookrightarrow N$ such that $S^n$ is homology trivial in $H^n (N,\mz_2)$ and the normal bundle of $S^n$ is nontrivial but stable trivial.
\end{coro}

\begin{coro}
If $H_n (M^{2n+2},\mz_2)=0,\ n$ odd, then there exists a quadratic form $\psi: H_{n}(N,\mz_2)\rightarrow \mz_2$ and for any embedding $\varphi: S^n \hookrightarrow N$ with $\varphi$ nullhomotopic in $M$, $\psi(\varphi_* ([S^n]))=0$ if and only if the normal bundle of $\varphi(S^n)$ is trivial.
\end{coro}

\begin{exmp}
For a smooth hypersurface $F_d \hookrightarrow \mc P^{n+1},\ n$ odd, $H^*(\mc P^{n+1},\mz_2)\cong\mz_2 [x]/(x^{n+2})$, deg$x$=2. $u_{F_d}=dx$. Then, if $d$ is even, $u_{F_d}=0$ and for the generator $x^{\frac{n+1}{2}} $ of $H^{n+1}(\mc P^{n+1},\mz_2)$:
\[
<Sq(x^{\frac{n+1}{2}}),[\mc P^{n+1}]>=<x^{\frac{n+1}{2}}(1+x)^{\frac{n+1}{2}},[\mc P^{n+1}]> = 1
\]
If $d$ is odd, $u_{F_d}=x$ and
\[
<\frac{Sq(x^{\frac{n+1}{2}})}{1+x},[\mc P^{n+1}]=<x^{\frac{n+1}{2}}(1+x)^{\frac{n-1}{2}},[\mc P^{n+1}]>=0
\]
We see in this case, the obstruction is just determined by the degree $d$ of $F_d$.
\end{exmp}

\subsection{topological blow up case}
We continue to discuss the topological blow up case. We still let $M^{2n+2}$ be a 1-connected, closed, smooth manifold and $i_N :N^{2n}\hookrightarrow M^{2n+2}$ be a closed submanifold of $M^{2n+2}$ with $\pi_i (M,N)=0,\ i<n+1.$

Let $i_Y:Y\hookrightarrow N$ be a closed submanifold of $N$ and suppose we can blow up $N$ along $Y$, we can also blow up $M$ along $Y$ since the normal bundle of $i_N (N)$ is a complex line bundle. Denote $\tilde{M}$ by the topological blow up of $M$ along $Y$ and $\tilde{N}$ by the topological blow up of $N$ along $Y.$
We obtain an embedding:
\[
\begin{CD}
\tilde{N} @>\widetilde{i_N}>> \tilde{M}\\
@VV \pi_N V       @VV\pi_MV \\
N @>i_N>>  M
\end{CD}
\]
We see $\tilde{N}$ is a submanifold of $\tilde{M}$ with codimension 2 and $\pi_1(\tilde{N})=\pi_1(\tilde{M})=1$, since dim$N-$dim$Y>2$. By excision theorem, $H_i (\tilde{M},\tilde{N},\mz)=H_i(M-Y,N-Y,\mz)=H_i (M,N,\mz)=0, \ i<n+1$. By relative Hurewicz theorem, $\pi_i (\tilde{M},\tilde{N})=0, \ i<n+1$. Thus we get:
\begin{lem}
$(\tilde{M},\tilde{N})$ is n-connected.
\end{lem}

In the sense of subsection 4.1, let $\tilde{L}$ be the complex line bundle with Euler class $e(\tilde{L})=PD^{-1}((\widetilde{i_N})_*[\tilde{N}])$ and $\tilde{W}=D(-\tilde{L})$,
we can talk about the obstruction $v_{n+1}(\tilde{W})\in H^{n+1}(\tilde{M},\mz_2)$. We need to determine $u_{\tilde{N}}$ first.

\begin{lem}
For $H^2(\tilde{M},\mz_2)=\pi_{M}^{*}H^2(M,\mz_2)\oplus \mz_2 \alpha_{j_{P_M}}^{*}(H^0(Y,\mz_2))$:
\[
u_{\tilde{N}}= \pi_{M}^{*}u_N + \alpha_{j_{P_M}}^{*}(1)
\]
\end{lem}
\begin{proof}
First, we know $\pi_M$ and $\pi_N$ are degree one maps and by the diagram:
\[
\begin{CD}
\tilde{N} @>\widetilde{i_N}>> \tilde{M}\\
@VV \pi_N V       @VV\pi_MV \\
N @>i_N>>  M
\end{CD}
\]
we see $(\pi_M)_* (\widetilde{i_N})_* [\tilde{N}]= (i_N)_* [N]$, $\alpha_{j_{P_{M}}}^* {u_{\tilde{N}}}=PD^{-1}((\pi_M)_*(\widetilde{i_N})_* [\tilde{N}] )=PD^{-1}(i_N)_* [N]=u_N$,
then we have $u_{\tilde{N}}= \pi_{M}^{*} u_N + a \alpha_{j_{P_M}}^{*}(1), \ a\in \mz_2$.

Second, by the diagram:
\[
\begin{CD}
\tilde{N} @>\widetilde{i_N}>> \tilde{M}\\
@Aj_{P_N}AA               @Aj_{P_M}AA \\
P(\eta_{Y}^{N})@>>>       P(\eta_{Y}^{M})
\end{CD}
\]
where $\eta_{Y}^{N}$ and $\eta_{Y}^{M}$ are normal bundles of $Y$ in $N$ and $M$.
We find $\tilde{N}$ intersects $P(\eta_{Y}^{M})$ transversally at $P(\eta_{Y}^{N})$ and $j_{P_N}^{*} \tilde{L}$ is isomorphic to the normal bundle of $P(\eta_{Y}^{N})\hookrightarrow P(\eta_{Y}^{M})$, whose mod 2 Euler class is, by calculation, just equal to $\overline{\xi_N} + i_Y^{*}u_N\in H^2(P(\eta_{Y}^{N}),\mz_2)=\mz_2\overline{\xi_N}\oplus H^2(Y,\mz_2)$. By the mod 2 Euler class of $j_{P_N}^{*} \tilde{L}$ is also equal to $j_{P_N}^{*}u_{\tilde{N}}=a\overline{\xi_N} + j_{P_N}^* \pi_M^*u_N= a\overline{\xi_N}+ i_{Y}^{*}u_N$. We see $a=1$ and we finally obtain:
\[
u_{\tilde{N}}= \pi_{M}^{*}u_N + \alpha_{j_{P_M}}^{*}(1)
\]
\end{proof}
So we get directly:
\begin{prop}
$v_{n+1}(\tilde{W})=0 \Longleftrightarrow$ for all $x\in H^{n+1}(\tilde{M},\mz_2)$
\[
<\frac{Sqx}{1+\pi_{M}^{*}u_N+ \alpha_{j_{P_M}}^{*}(1)},[\tilde{M}]> =0
\]
\end{prop}

\begin{exmp}
We consider a simple but common example. Let $M^{2n+2}$ be a 1-connected, closed, smooth manifold and $i_N :N^{2n}\hookrightarrow M^{2n+2}$ be a closed submanifold of $M^{2n+2}$ with $\pi_i (M,N)=0,\ i<n+1.$ Choose one point $p\in N$, we blow up $M$ and $N$ at this point to obtain $\tilde{M}$, $\tilde{N}$ and the embedding $\tilde{N}\longrightarrow \tilde{M}$. For two obstructions $v_{n+1}(W)$ and $v_{n+1}(\tilde{W})$, we claim:
\[
v_{n+1}(W)=0 \Longleftrightarrow v_{n+1}(\tilde{W})=0
\]
Indeed, we see $H^* (\tilde{M})=\pi_{M}^{*}H^*(M) \oplus \alpha_{j_{P_M}}^{*}(H^* (\mc P^n ))$. The cohomology group $H^* (\mc P^{n})\cong \mz_2 [x]/(x^{n+1})$ and we get:
\[
H^{n+1}(\tilde{M})=\pi_{M}^{*}H^{n+1}(M)\oplus \mz_2 \alpha_{j_{P_M}}^{*}(x^{\frac{n-1}{2}})
\]
Let $y:=\alpha_{j_{P_M}}^{*}(1)\in H^2(\tilde{M})$, we have $y^{\frac{n+1}{2}}=\alpha_{j_{P_M}}^{*}(x^{\frac{n-1}{2}})$ and $y\cup \pi_{M}^{*}a=0,\ 1\neq a\in H^{*}(M)$.

For $y^{\frac{n+1}{2}},$ we see
\[
\frac{1}{1+\pi_{M}^{*}u_N +y}=\frac{1}{1+y} + \frac{\pi_{M}^{*}u_N}{(1+y)(1+\pi_{M}^{*}u_N+y)}
\]
we have:
\begin{equation*}
\begin{split}
  & <\frac{Sqy^{\frac{n+1}{2}}}{1+\pi_{M}^{*}u_N+y},[\tilde{M}]>\\
= & <\frac{Sqy^{\frac{n+1}{2}}}{1+y},[\tilde{M}]>+<\frac{\pi_{M}^{*}u_N Sqy^{\frac{n+1}{2}}}{(1+\pi_{M}^{*}u_N+y)(1+y)},[\tilde{M}]>\\
= & 0
\end{split}
\end{equation*}
here $y\pi_{M}^{*}u_N =0$ implies $\pi_{M}^{*}(u_N) Sqy^{\frac{n+1}{2}}=0$.

For any $a\in H^{n+1}(M)$, we have:
\begin{equation*}
\begin{split}
  & <\frac{\pi_{M}^{*}(a)}{1+\pi_{M}^{*}u_N+y},[\tilde{M}]>\\
= & <\frac{\pi_{M}^{*}(a)}{1+\pi_{M}^{*}u_N},[\tilde{M}]>+<\frac{y\pi_{M}^{*}(a)}{(1+\pi_{M}^{*}u_N+y)(1+y)},[\tilde{M}]>\\
= & <\frac{a}{1+u_N},[M]>
\end{split}
\end{equation*}
Thus we get our claim.

\end{exmp}

%{\bf Acknowlegement} I would like to thank...

\vspace{1cm}
\noindent{Wang Wei}\\
School of Mathematical Sciences,\\
Fudan University,\\
Shanghai 200433, P.R. China.\\
Email: weiwang@amss.ac.cn


\begin{thebibliography}{50}
\setlength{\itemsep}{-3pt}
\small
\bibitem{Browder} Browder, William Surgery on simply-connected manifolds. Ergebnisse der Mathematik und ihrer Grenzgebiete, Band 65. Springer-Verlag, New York-Heidelberg, 1972. ix+132 pp
\bibitem{Browder2} Browder, William Complete intersections and the Kervaire invariant. Algebraic topology, Aarhus 1978 (Proc. Sympos., Univ. Aarhus, Aarhus, 1978), pp. 88¨C108,
Lecture Notes in Math., 763, Springer, Berlin, 1979.
\bibitem{Duan} Haibao Duan, Banghe Li. Topology of Blow-ups and Enumerative Geometry. arXiv:0906.4152v4 [math.AG]
\bibitem{Geiges} Geiges, Hansjorg; Pasquotto, Federica A formula for the Chern classes of symplectic blow-ups.  J. Lond. Math. Soc. (2)  76  (2007),  no. 2, 313¨C330.
\bibitem{Griff}Griffiths, Phillip; Harris, Joseph Principles of algebraic geometry. Reprint of the 1978 original. Wiley Classics Library. John Wiley Sons, Inc., New York, 1994.
\bibitem{Lascu}   Lascu, A. T.; Scott, D. B. A simple proof of the formula for the blowing up of Chern classes.  Amer. J. Math.  100  (1978), no. 2, 293¨C301.
\bibitem{Macduff} McDuff, Dusa Examples of simply-connected symplectic non-K0Š1hlerian manifolds.  J. Differential Geom.  20  (1984),  no. 1, 267¨C277.
\bibitem{Milnor} Milnor, John W.; Stasheff, James D. Characteristic classes. Annals of Mathematics Studies, No. 76. Princeton University Press, Princeton, N. J.; University of Tokyo Press, Tokyo, 1974. vii+331 pp.
\bibitem{Voisin} Voisin, Claire Hodge theory and complex algebraic geometry. I. Translated from the French by Leila Schneps. Reprint of the 2002 English edition. Cambridge Studies in Advanced Mathematics, 76. Cambridge University Press, Cambridge, 2007.

\end{thebibliography}
\end{document}